
\magnification=1200
\def\R{\hbox{$I\kern-3.5pt R$}}
\def\N{\hbox{$I\kern-3.5pt N$}}
\def\V{\Vert}
\font\fp=cmr8
\def\P{{\noindent \it Proof.- }}
\def\fin{ Q.E.D. \medskip}

\def\var{\overline\va}

\hsize=15true cm
\hoffset=0.5true cm
\vsize=23true cm
\def\va{\varphi}
\centerline{\bf INTERPOLATION OF OPERATORS WHEN THE EXTREME}
\centerline {\bf SPACES ARE $L^\infty$}
\bigskip
\centerline {by}
\bigskip
\centerline{ Jes\'us Bastero\footnote
*{\fp Research partially supported by DGICYT PS87-0059} and Francisco J.
Ruiz\footnote {**}{\fp Research partially supported by DGICYT
PB89-0181-C02-02} }
\bigskip
\midinsert
\narrower\narrower
\noindent  ABSTRACT. {\sl In this paper, equivalence between interpolation
properties of linear operators and monotonicity conditions are studied, for a
pair $(X_0,X_1)$ of rearrangement invariant quasi Banach spaces, when the
extreme spaces of the interpolation are $L^\infty$ and a pair $(A_0,A_1)$
 under
some assumptions. Weak and restricted weak intermediate spaces fall in our
context. Applications to classical Lorentz and Lorentz-Orlicz spaces are
 given.}
\endinsert
\bigskip

\beginsection \$ 0. Introduction.

 Let $(A_0,A_1)$, $(B_0,B_1)$, $(X_0,X_1)$ be three pairs of rearrangement
invariant Banach function spaces (see definitions below) over the interval
 $I$
($I=[0,1]$ or $[0,\infty )$).
 Let ${\cal A}((A_0,A_1),(B_0,B_1))$ denote the class of linear (or
 quasilinear
or Lipschitz) operators which are bounded from $A_0$ into $A_1$ and from
 $B_0$
into $B_1$.  The pair $(X_0,X_1)$ is said to have  the {\it linear (or
quasilinear or Lipschitz) interpolation property with respect to the class
${\cal A}$} if every member of ${\cal A}$ can be extended to a bounded
operator
from $X_0$ into $X_1$.

This interpolation property has been extensively studied in its
 connection with
many aspects concerning  r.i. spaces, for instance, Boyd or Zippin's
indexes, monotonicity conditions, boundedness of some suitable \lq\lq
maximal"
operators and so on. Here we are  concerned with the case
$B_0=B_1=L^\infty$ and particularly in  connection with the
 monotonicity property $(\cal M)$ given in \$ 1 and the boundedness of
only one
operator.

 In this direction the former result is contained in Calderon's
paper {\bf [5]} where it is showed  that  both properties, say linear
interpolation and monotonicity, are equivalent in the  case of $A_0=A_1=L^1$.
Later on,  Lorentz and Shimogaki {\bf [10]} extended this result to the case
$A_0=A_1=L^p$ with $p>1$. The technique used by them consists on a
linearization  process of the $L^p$ case.

Sharpley, Maligranda and  other autors  (see {\bf [11]} and references
quoted there) studied the case $A_0=\Lambda (X)$, $A_1=M(X)$ (see definitions
 in \$ 2)
and  $B_0 = B_1 = L^\infty$ or $B_0=\Lambda (Y)$,  $B_1=M (Y)$ relating the
interpolation properties with the boundedness of  only one \lq\lq maximal"
operator ({\bf [18, theorem 4.7]}, {\bf [11, theorem 4.5]}). On the other
 hand,
Maligranda {\bf[11]} obtained equivalence between the interpolation property
 for
Lipschitz operators and monotonicity condition in the case $A_0=\Lambda (X)$,
$A_1=M(X)$ and $B_0=B_1=L^\infty$. When $X=L^p$, $p>1$,  then $\Lambda
(X)=L^{p,1}$ and $M(X)= L^{p,\infty}$. So we can see Maligranda's result is
close to Lorentz-Shimogaki's one. The spaces with the interpolation property,
when the extreme spaces are $\Lambda (X)$ and $M(X)$, are generally known in
 the
literature as {\it weak type intermediate spaces}.

These papers leave out the more \lq\lq natural" case where
$A_0=L^p$, $A_1=L^{p,\infty}$ or, more generally, $A_0=X$, $A_1=M(X)$.
In fact,
following the usual terminology in Fourier Analysis, it should be reserved
the
term   {\it weak type intermediate spaces} to  spaces having  the
interpolation
property in this last setting, while the spaces with the interpolation
property
in the setting before stated should be named {\it restricted weak type
intermediate spaces}.

Our final purpose is to study this  \lq\lq intermediate" case between
Lorentz-Shimo\-ga\-ki's and Maligranda's. In order to do that, our main
 tool
consists of obtaining, in a very general context, equivalence between
interpolation properties of linear, quasilinear or Lipschitz type and
monotonicity condition $(\cal M)$. When this result  is established it is an
 easy
consequence to reduce the linear interpolation property to the boundedness of
only one quasilinear operator.

This general result can be applied in the both cases   stated before, namely,
weak and restricted weak intermediate. So, on one hand we obtain some
generalizations of Maligranda's results and on the other one we obtain several
results in the case of $A_0=X$, $A_1=M(X)$. When $A_0=L^p$, the quasilinear
operator can be iterated and, as a consequence, we obtain that the weak type
intermediate spaces are exactly the restricted weak intermediate spaces.

Moreover, by using a characterization  about the boundedness of the Hardy
operator in Lorentz spaces due to Ari\~no and Muckenhoupt, we can characterize
the  Lorentz spaces which are intermediate in terms of handly conditions on
the
weights. Finally, the last part of the paper is devoted to extend some of the
previous results to the more general case of Lorentz-Orlicz spaces.

The paper is organized in two sections: the first one
contains the notations and the general results and  the second one the
applications.
 \bigskip
\beginsection \$ 1. General results

A Banach space $(X,\Vert . \Vert)$ of real-valued, locally integrable,
Lebesgue
measurable functions on $I$ ($I=[0,1]$ or $[0,\infty )$) is said
to be a rearrangement invariant Banach function space over
$I$ (in short r.i. space) if it satisfies the following conditions:

\item{\bf i)} If $|g|\leq |f|$ a.e. and $f\in X$, then $g\in X$ and
$\Vert g\Vert \leq \Vert f\Vert$.

\item{\bf ii)} $0\leq f_n \uparrow$, $\sup_{n\in \N}\Vert f_n\Vert\leq M$,
 imply
that $f=\sup f_n \in X$ and $\Vert f\Vert  =\sup_{n\in \N}\Vert f_n\Vert$.

\item{\bf iii)} $X$ contains the simple integrable functions.

\item{\bf iv)} $f\in X \Longleftrightarrow f^\ast\in X$ and $\Vert f \V = \V
f^\ast \V$, where $f^\ast$ denotes the nonincreasing rearrangement of the
function $f$. \medskip

Fact ii) is known in the literature as Fatou property (cf. {\bf [9]}). It is
quite clear that if $X$ is r.i. then  $L^1\cap L^\infty \hookrightarrow
X\hookrightarrow L^1+L^\infty$ (where the symbol $\hookrightarrow$ signifies
continuously embedded).

 A classical result by Lorentz and
Luxemburg   ensures that for these spaces
\item{\bf v)} $\displaystyle \V f\V = \sup_{\V g\V_{X'}\leq 1} \left | \int_I
 fg
\right |$, where  $X'$ is the associated
space of $X$ which is also r.i. space. In particular $X=X''$ isometrically.
\medskip

The fundamental function $\phi_X$ of an r.i. space is defined by
$$\phi_X(t)=\V\chi_{[0,t]}\V,  t\in I.$$

There is no loss of generality if we assume for  $\phi_X$ to be
positive, nondecreasing, absolutely continuous far from the origin, concave
and to verify (see {\bf [18]}, {\bf [21]}):

\item{\bf vi)} $\phi_X(t)\phi_{X'}(t)=t$, for all $t\in I;$

\item{\bf vii)} $\displaystyle {{d\phi_X(t)}\over {dt}} \leq {{\phi_X(t)}\over
t}$, a.e. on $I$ .

In which follows it may be convenient to let $X$  be a quasi-Banach r.i.. The
main difference occurs in the triangle inequality satisfied in $X$, i.e.
$\V f+g\V \leq C(\V f\V +\V g\V )$, for some constant $C\geq 1$. In this case
we suppose that a quasi-Banach  space $X$ satisfies properties i), ii), iii),
iv) but, in general, no other conditions  will be assumed.
We say that a quasi-Banach function space is $\sigma$-order continuous
 if   every order bounded nondecreasing sequence converges in the quasi-norm
topology (cf. {\bf [9, Proposition 1.a.8]}).

\medskip
 Let $(A_0,A_1)$, $(X_0,X_1)$, $i=1, 2,$ be two couples of r.i. quasi-Banach
spaces on $I$ such
 that $A_i\cap L^\infty \hookrightarrow X_i\hookrightarrow A_i+L^\infty$, $i=1,
2.$
 We say that the couple $(X_0,X_1)$ belongs to:

- ${\cal LI}(A_0,A_1;L^\infty )$  if any  linear operator $T:A_0+L^\infty
 \rightarrow A_1+L^\infty$ which is bounded from $A_0$ into $A_1$ and from
$L^\infty$ into $L^\infty$ is also bounded from $X_0$ into $X_1.$ The closed
graph theorem implies that there exists a constant $C\geq 1$ such that
$$
\V T\V_{X_0\rightarrow X_1} \leq C\max\{\V T\V_{A_0\rightarrow A_1},
\V T\V_{L^\infty\rightarrow L^\infty}\}.
$$

- ${\cal QLI}(A_0,A_1;L^\infty )$  if any  quasilinear operator
 $T:A_0+L^\infty
 \rightarrow A_1+L^\infty$ which is bounded from $A_0$ into $A_1$ and from
$L^\infty$ into $L^\infty$ is also bounded from $X_0$ into $X_1.$

- ${\cal LPI}(A_0,A_1;L^\infty )$  if any  operator $T:A_0+L^\infty
 \rightarrow A_1+L^\infty$ that is a  Lipschitz operator from
$A_0$ into $ A_1$ and from $L^\infty$ into
$L^\infty$  also maps $X_0$ into $X_1.$

Recall that a map $T:A_0\rightarrow A_1$ is  bounded quasilinear if there are
constants $K, C\geq 1$ such that $|T(\lambda f)|=|\lambda ||T(f)|$,
 $|T(f+g)|\leq
K(|T(f)|+|T(g)|)$ and $\V T(f)\V \leq C\V f\V$. We define
 $\V T\V_{A_0\rightarrow A_1} = \inf C$.
In the same way, a map \ \  $T:A_0\rightarrow A_1$ is a  Lipschitz operator if
there is a constant $C\geq 1$ such that $T(0)=0$, $\V Tf-Tg\V\leq C\V f-g\V $;
 we
define $\V T\V_{A_0\rightarrow A_1} = \inf C$.

Next we shall introduce another class of spaces and for that we need the
following
\medskip
\proclaim{Lemma 1}. If $X$ is a quasi-Banach r.i. space and $f\in X+L^\infty$
and $m(E)<\infty$ then $f\chi_E\in X$.

\P  Let $f=g+h$ with $g\in X$ and $h\in L^\infty$. Since $g\chi_E\in X$ and
$|h\chi_E|\leq \V h\V_\infty \chi_E \in X$ the result follows immediately.
\fin

We say that the couple $(X_0,X_1)$ belongs to:

- ${\cal M}(A_0,A_1)$  if there exists a constant $C\geq 1$ such that
if $f\in X_0$, $g\in A_1+L^\infty$ and
$$
\V f^\ast\chi_{[0,t]}\V_{A_0}\geq
\V g^\ast\chi_{[0,t]}\V_{A_1} \ \ \ \forall t > 0\eqno {(\cal M)}$$
then $g \in X_1$ and $\V g\V_{X_1}\leq C\V f\V_{X_0}$.
\medskip
It is clear that
$$
{\cal QLI}(A_0,A_1;L^\infty ) \cup
 {\cal LPI}(A_0,A_1;L^\infty )\subseteq {\cal LI}(A_0,A_1;L^\infty ).
$$
Under some more restrictive assumptions the four classes of maps
introduced before coincide.
  \medskip
\proclaim {Proposition 1}. Let $(A_0,A_1)$ a couple of quasi-Banach r.i.
spaces such that $\phi_{A_1}(t)\leq C\phi_{A_0}(t)$ for all $t>0$, then
$$
\leqno {\bf (1.1)} \ {\cal M}(A_0,A_1)\subseteq {\cal QLI}(A_0,A_1;L^\infty )
$$
{\bf (1.2)} If $I=[0,1]$ and $A_0$ is $\sigma$-order continuous then
$${\cal M}(A_0,A_1)\subseteq
{\cal LPI}(A_0,A_1;L^\infty ).$$

\P {\bf (1.1)} Let $(X_0,X_1)$ be a couple in ${\cal M}(A_0,A_1)$
and let $T$ be a quasilinear map  $T:A_0+L^\infty \rightarrow
 A_1+L^\infty$, bounded from $A_0$ into $A_1$ and from $L^\infty$ into
$L^\infty$. Suppose that $\V T\V_{A_0\rightarrow A_1} \leq 1$ and
$\V T\V_{L^\infty\rightarrow L^\infty} \leq 1$. We have to show that $T$ is
bounded from $X_0$ into $X_1$. In order to do it
let $f$ be an element in $X_0$.
We only need to  prove that
$$
\V (Tf)^\ast\chi_{[0,t]}\V_{A_1}\leq \V (Cf)^\ast
\chi_{[0,t]}\V_{A_0} \eqno
{(\ast)} $$
for all $t>0$.

 We know that
$$
\V (Tf)^\ast\chi_{[0,t]}\V_{A_1} = \sup\V (Tf)\chi_E\V_{A_1}
$$
where $E$ runs over the borelians in $I$ with $m(E)\leq t$.  Set
$s=f^\ast (t)$
and define
$$f_{(s)}=\cases {
s, & if $f>s$ ;\cr
-s, & if $f<-s$ ; \cr
f, & otherwise  \cr }$$
and $f^{(s)}=f-f_{(s)}$. Since $f_{(s)}\in L^\infty$  and $f^{(s)}\in A_0$ we
have
 $$\eqalign{
\V (Tf_{(s)})\chi_E \V_{A_1} &\leq \V Tf_{(s)}\V _{\infty} \V \chi_E\V_{A_1}
 \cr
& \leq \V f_{(s)}\V_\infty \phi_{A_1}(t)\leq Cs\phi_{A_0}(t) \cr
&\leq C\V f^\ast (t)\chi_{[0,t]}\V_{A_0}\leq C\V f^\ast \chi_{[0,t]}\V_{A_0}.
}$$
Now
$$\eqalign{
\V (Tf^{(s)})\chi_E \V_{A_1} &\leq \V Tf^{(s)}\V _{A_1} \cr
& C\leq \V f^{(s)}\V_{A_0} \leq C \V f^\ast \chi_{[0,t]}\V_{A_0}.
}$$
and hence we easily obtain the inequality ($\ast$).

\noindent {\bf (1.2)} In order to show now  that the couple $(X_0,X_1)$
 belongs
to ${\cal LPI}(A_0,A_1;L^\infty )$ we will follow the ideas of {\bf [11,
 lemma
4.4]}. For the sake of completeness, we include the proof here. Suppose
that
$T$  is a Lipschitz operator mapping $A_0+L^\infty$ into $A_1+L^\infty$
with $\V T\V_{A_0\rightarrow A_1} \leq 1$ and  $\V T\V_{L^\infty\rightarrow
L^\infty}  \leq 1$.

If $0<t\leq 1$ and $f\in X_0$, we set $f^\ast (t)=s$. Consider  $f_{(s)}$ and
$[(Tf)^\ast]_{(s)}$  defined in a similar way as before. Since $f_{(s)}\in
L^\infty$ we have $\V T(f_{(s)})\V_\infty\leq s$ and so,
$$ |Tf-(Tf)_{(s)}|\leq |Tf-T(f_{(s)})|.$$
 If $0<x\leq t$,
$$|(Tf)^\ast (x)|\leq |(Tf)^\ast (x)-[(Tf)^\ast]_{(s)}(x)| +
[(Tf)^\ast]_{(s)}(x)$$
thus
$$\eqalign{
\left\V (Tf)^\ast \chi_{[0,t]}\right\V_{A_1} &\leq C\left ( \left\V \left
|(Tf)^\ast -[(Tf)^\ast ]_{(s)}\right |\chi_{[0,t]}\right\V_{A_1} +
 f^\ast (t)\V \chi_{[0,t]}\V_{A_1} \right)
 \cr
&\leq C\left ( \left\V \left[ Tf
-(Tf)_{(s)}\right ]\chi_{[0,t]}\right\V_{A_1} +
 \V f^\ast  \chi_{[0,t]}\V_{A_0} \right)\cr
&\leq C\left ( \left\V |Tf
-T(f_{(s)})|\chi_{[0,t]}\right\V_{A_1} +
 \V f^\ast  \chi_{[0,t]}\V_{A_0} \right )
\cr
&\leq C\left ( \left\V Tf
-T(f_{(s)})\right\V_{A_1} +
 \V f^\ast  \chi_{[0,t]}\V_{A_0} \right )
\cr
&\leq
C\left ( \V (f
-f_{(s)})\V_{A_0} +
 \V f^\ast  \chi_{[0,t]}\V_{A_0} \right )\cr
&\leq
C
 \V f^\ast  \chi_{[0,t]}\V_{A_0}
}$$
 This implies  that $T$ maps  $X_0$ into $X_1$. Next by using the fact that
$A_0$ is $\sigma$-order continuous    we can follow the proof of the  theorem
4.5 from {\bf [11, theorem 4.5]}, and conclude the proof of this part. (The
constants $C$ appearing above may change from line to line). \fin
\medskip

In order to obtain another implications we need to restric our attention to
Banach spaces.
\medskip

\proclaim Proposition 2. Let $A_0$ be a Banach r.i. space and
let $A_1$ be a quasi-Banach r.i. space. Suppose that
$\phi_{A_0}\leq C\phi_{A_1}$ and $\displaystyle {1\over
{\phi_{A_1}}}\in A_1$. Then
 $${\cal LI}(A_0,A_1;L^\infty
)\subseteq {\cal M}(A_0,A_1).$$

\P Let $f\in X_0$, $g \in A_1+L^\infty$ such that
$$
\V f^\ast\chi_{[0,t]}\V_{A_0}\geq
\V g^\ast\chi_{[0,t]}\V_{A_1}, \ \ \forall t>0
$$
For every $k\in {\bf Z}$ we define $E_k=[2^k,2^{k+1})$. We have
$$\eqalign{
g^\ast (t) & =\sum^\infty_{-\infty} g^\ast (t)\chi_{E_k}(t) \cr
&\leq \sum^\infty_{-\infty} g^\ast (2^k)\chi_{E_k}(t) \cr
&\leq  \sum^\infty_{-\infty} {1\over {\phi_{A_1}(2^k)}}\
\V f^\ast    \chi_{[0, 2^k]}\V_{A_0} \   \chi_{E_k}(t)
}$$
For every $k\in {\bf Z}$ we can choose a function $h_k\in A_0'$ with $\V
h_k\V_{A'_0}\leq 1$, such that
  $$\V f^\ast\chi_{[0,2^k]}\V_{A_0}\leq 2\int_I
f^\ast\chi_{[0,2^k]} h_k.$$
Then
$$g^\ast (t) \leq 2 \sum^\infty_{-\infty} {1\over {\phi_{A_1}(2^k)}}
\left(\int_0^{2^k}
f^\ast h_k\right)
\chi_{E_k}(t).$$
For any locally integrable function $\varphi$  on $I$ we define the
\lq\lq linear" operator $T$ by
 $$T\varphi (t)= 2\sum^\infty_{-\infty} {1\over {\phi_{A_1}(2^k)}}
\left(\int_0^{2^k}
\varphi h_k\right)
\chi_{E_k}(t).$$
It is clear that if $\varphi \in L^\infty$ and $t\in E_k$ we have
$$\eqalign{
|T\varphi (t)|&\leq  {2\over {\phi_{A_1}(2^k)}}\V\varphi\V_\infty
\V
\chi_{[0,2^k]}\V_{A_0}\cr
&\leq C\V\varphi\V_\infty.}$$

On the other hand, if $\varphi\in A_0$ then
$$|T\varphi (t)|\leq 2\V\varphi\V_{A_0} \sum^\infty_{-\infty} {1\over
{\phi_{A_1}(2^k)}}
\chi_{E_k}(t).$$
For $t\in E_k$ the triangle inequality of the quasi-norm implies that
$$
\phi_{A_1}(t)\leq \phi_{A_1}(2^{k+1})\leq C \phi_{A_1}(2^k)$$
hence
$$\eqalign{
|T\varphi (t)|&\leq C\V\varphi\V_{A_0}\sum_{-\infty}^\infty{
{\chi_{E_k}(t)}\over{
\phi_{A_1}(t)}}  \cr
&=C\V\varphi\V_{A_0}{1\over {
\phi_{A_1}(t)}}
}. $$
Since
$\displaystyle {1\over
{\phi_{A_1}}}\in A_1$ we obtain that $T$ is also bounded from $A_0$ into $A_1$.
 Eventually  we have that  $g\in A_1$ and $\V g\V_{A_1}\leq C\V  f \V_{A_0}$
because  $g^\ast\leq Tf^\ast\in A_1.$
\fin
\medskip
Among the  results above we want to emphasize the following

\proclaim Theorem 1. Let $A_0$ be a Banach r.i. space and let $A_1$ be a
quasi-Banach r.i. space. Suppose that $C^{-1}\phi_{A_1}\leq \phi_{A_0}\leq
C\phi_{A_1}$ for some constant $C$ and $\displaystyle {1\over {\phi_{A_1}}}
\in
A_1$. Then
$$
{\cal LI}(A_0,A_1;L^\infty)= {\cal M}(A_0,A_1).
$$
Moreover, a couple of quasi Banach r.i. spaces $(X_0,X_1)$ belongs  to any of
this classes if and only if the quasilinear operator
$$
Q\varphi (t)= { 1\over {\phi_{A_1}(t)}} \V \varphi\chi_{[0,t]}\V_{A_0}
$$
 is bounded from $X_0$ into $X_1$ for nonincreasing functions.

\P The first part follows from propositions 1 and 2. For the second part,
observe that $Q$  is bounded from $A_0$ into $A_1$ and from $L^\infty$ into
$L^\infty$ and so if $(X_0,X_1) \in {\cal LI}(A_0,A_1;L^\infty) (=
{\cal QLI}(A_0,A_1;L^\infty)$, then $Q$ is bounded from $X_0$ into $X_1$. On
the other hand, note that condition $(\cal M)$ implies that $g^\ast
(t)\leq Qf^\ast (t)$ and, therefore, the boundedness of $Q$ for nonincreasing
functions implies that $(X_0,X_1) \in {\cal M}(A_0,A_1;L^\infty)$\fin
\medskip
{\noindent \bf Remarks.}
\noindent {\bf i)}  Under the same hypotheses as in
proposition 2, but supposing only that $A_0$ is  quasi-Banach, we can prove in
a simpler way that  ${\cal QLI}(A_0,A_1;L^\infty ) \subseteq{\cal M}(A_0,A_1)$.
The operator we have to use instead of $T$  is $Q$.
\medskip
\noindent {\bf ii)}  Proposition 2 is not true when $A_0$ is
quasi-Banach. For instance, let $I=[0,1]$, $A_0=A_1=L^{p,\infty}$, $0<p\leq 1$.
In this case it is easy to see that the couple $(L^p,L^p) \notin
{\cal M}(L^{p,\infty},L^{p,\infty};L^{\infty})= {\cal
QLI}(L^{p,\infty},L^{p,\infty};L^{\infty})$, but by using a result by Kalton
 {\bf
[7, Theorem 1.1]}, we can deduce that $(L^p,L^p)\in
 {\cal LI}(L^{p,\infty},L^{p,\infty};L^{\infty}),$ (this result was quoted to
the authors by Oscar Blasco).
 \bigskip

\beginsection \$ 2. Applications.

Not all Banach r.i. spaces $A_1$ satisfy the condition $\displaystyle {1\over
{\phi_{A_1}}}\in A_1.$ In order to study this property we
introduce the Lorentz spaces as they appear in {\bf [2]}, {\bf [11]}, {\bf
[18]}, {\bf [21]}.

 In what follows we assume that $X$ is a Banach r.i. space. We denote by:
\medskip
- $\Lambda (X)$  the space of all measurable functions  with
$$
\V f\V_{\Lambda (X)}=\int_If^\ast (t)d\phi_X (t)<\infty.
$$
Since $\phi_X$ is concave, the expression
$\V f\V_{\Lambda (X)}$ is a norm and moreover  $\Lambda (X)$ is a Banach r.i.
space.
\medskip
- $M(X)$  the space of all measurable functions $f$ for which there exists
$f^{\ast\ast}$ and
$$
\displaystyle \V f\V_{M(X)}=\sup_{t\in
I}\phi_X(t)f^{\ast\ast} (t)<\infty.
$$
Recall that $f^{\ast\ast}$, the {\it Hardy transform} of $f^\ast$, is defined
 by
$$H(f^\ast)(t)= f^{\ast\ast}(t)={1\over t} \int_0^t f^\ast.
$$
$M(X)$ is again a Banach r.i. space.
\medskip
-  $M^\ast (X)$  the space of all measurable functions for which
$$
 \V f\V_{M^\ast (X)}=\sup_{t\in I}\phi_X(t)f^\ast (t)<\infty.
$$
The function $\V .\V_{M^\ast (X)}$ is a quasinorm on $M^\ast (X)$.

It is clear that for these spaces we have:
\medskip
 \noindent {\bf i)} $X\subseteq M^\ast (X)$
\ , \  $\Lambda (X)\hookrightarrow X \hookrightarrow M(X)$\ , \
$M(X) \subseteq M^\ast (X)$

\noindent {\bf ii)} $\phi_{\Lambda (X)}= \phi_{M^\ast (X)}=\phi_X=\phi_{M(X)}$,

\noindent {\bf iii)}  ${1\over \phi_X}\in M^\ast (X)$,

\noindent {\bf iv)} $ {1\over \phi_X}\in X \Longleftrightarrow X=M^\ast (X)$.

\proclaim Lemma 2. Let $X$ be a  Banach r. i. space. The  following
conditions are equivalent:

{\sl \noindent{\bf (a)} The space $M^\ast (X)$ is
convexifiable (i.e. there is a norm on $M^\ast (X)$ equivalent to $\V
\cdot\V_{M^\ast (X)}$)

\noindent{\bf (b)} $M(X)=M^\ast (X)$

\noindent{\bf (c)} $\displaystyle {1\over
{\phi_X}}\in M (X)$

\noindent{\bf (d)} There exist  a constant $C>0$ such that
$$ {\phi_X(t)\over t}\int_0^t{ds\over\phi_X(s)}\leq C, \ \ \forall t\in I$$

\noindent{\bf (e)} $\V f\V_{M(X)} \sim \V f\V_{M^\ast (X)}, \ \forall f\in
 M^\ast
(X)$.}

\medskip
 \P We only will sketch {\bf (a)} $\Rightarrow$ {\bf (b)}. We may assume
there is a r.i. norm $|||.|||$ on $M^\ast (X)$ equivalent to $\V .\V_{M^\ast
(X)}$.  If $f\in
M^\ast (X)$ and $t>0$
$$f^{\ast\ast}(t)\phi_X(t)\leq C|||f^{\ast\ast}(t)\chi_{[0,t]}|||\leq
C|||f^\ast\chi_{[0,t]}|||$$
and therefore
$\V  f \V_{M(X)}\leq C\V  f \V_{M^\ast(X)}$.\fin
\medskip
Now we can state the results of preceding section in the framework of Lorentz
spaces.

\medskip

\proclaim Proposition 3. If $X$ is a r.i. Banach space then ${\cal
LI}(\Lambda (X), M^\ast (X); L^\infty )={\cal
QLI}(\Lambda (X), M^\ast (X); L^\infty )= {\cal
M}(\Lambda (X),M^\ast (X))={\cal
LPI}(\Lambda (X),M^\ast (X);L^\infty )$ (the last fact if $I=[0,1]$).
Furthermore, a couple of spaces  $(Y,Z)$ belongs to any of the classes before
stated if and only if the quasilinear operator $Q_{\Lambda (X)}$ defined by
$$Q_{\Lambda (X)}\varphi (t)= { 1\over {\phi_X(t)}} \int_0^t
\varphi (x)d\phi_X (x) \eqno {(2.1)}$$
is bounded from $Y$ into $Z$ for nonincreasing functions.

{\noindent \bf Notes}. {\bf (i)} This result has been already obtained by
Maligranda {\bf [11]} in the case $Y=Z$ .
\medskip
\noindent {\bf (ii)} The operator appearing in Proposition 3 is actually
$$
Q_{\Lambda (X)}\varphi = H(\varphi \circ \phi_X^{-1})\circ \phi_X.
$$

\bigskip
The preceding proposition can be applied to the class of
classical Lorentz spaces  $\Lambda (w,q)$ with non monotone weights.
Let $w$ be an a.e.
positive weight defined on $I=[0,\infty )$ such that $\int_0^tw<\infty,
\forall t<\infty $ and $\int_0^\infty w=\infty$. We recall that the classical
Lorentz space $\Lambda (w,q)$, $0<q\leq\infty$ is the class of all real
valued measurable functions on $I$ such that
$$\V f\V_{\Lambda (w,q)}=\cases
{
\left(\int_If^\ast (t)^qw(t)dt\right)^{1/q}<\infty , &if $0< q<\infty$ ;\cr
\sup_{t>0} f^\ast (t)w(t)<\infty , & if $q=\infty .$ \cr
}$$
For $q=\infty$ we  will only consider nondecreasing weights $w$.  Ari\~no
and Muckenhoupt {\bf [1]} have showed that given $0<q<\infty$, there exists a
constant $C>0$ such that the Hardy operator verifies $$
\V Hf\V_{\Lambda (w,q)}\leq C\V f\V_{\Lambda (w,q)}
$$
for all nonnegative and nonincreasing functions $f$ on $\R$  if
and only if the weight $w$ satisfies
 $$
\int_t^\infty{w(x)\over x^q}dx\leq{B\over t^q}\int_0^tw(x)dx \eqno{(AM_q)}
$$
for some constant $B>0$ and for all $t>0.$ Moreover, for $1\leq
q<\infty$, condition $(AM_q)$ implies that $\Lambda (w,q)$ is a
Banach space. Sawyer  {\bf [17]} proved that the converse is true for
$1<q<\infty$. Raynaud gave also another equivalent condition to this last fact
by using quasi-concavity conditions for the function  $W(t)=\int_0^tw(x)dx$
(see
{\bf [16]}).

In the case $q=\infty$ and $w$ nondecreasing, the same reasons appearing in
Lemma 2 show that $\Lambda (w,\infty)$ is a Banach space if and only if the
weight $w$ satisfies
$$
{w(t)\over t}\int_0^t {dx \over w(x)}\leq C \eqno (A_1)$$
for some constant $C>0$ and for all $t>0.$

If we suppose that the weights satisfy the conditions $(AM_q)$ or $(A_1)$ then
 $$
\V f^{\ast\ast}\V_{\Lambda (w,q )}\sim  \V f\V_{\Lambda (w,q)}
$$
and reciprocally.

In next statements, when we say that $\Lambda (w,q)$ is a Banach space we will
mean that conditions $(AM_q)$ or $(A_1)$ are satisfied.

By using Ari\~no-Muckenhoupt's result  stated  before we are going to be able
to
characterize the Lorentz spaces which are interpolated between $\Lambda (X)$,
$M^\ast (X)$ and $L^\infty.$

\proclaim Proposition 4. Let $X$ be a r.i. Banach space and suppose that
$\Lambda (w,q)$ is a Banach space. Then the following assertions are true:
\item{\bf (i)} For $1\leq q <\infty$,  $\Lambda (w,q)\in {\cal
LI} (\Lambda (X), M^\ast (X); L^\infty )
$ if and only if
$$
\int_t^\infty{w(x)\over \phi_X(x)^q}dx \leq
{B\over \phi_X(t)^q}\int_0^tw(x)dx, \qquad \forall t>0 \eqno {(2.2)}
$$
\item{\bf (ii)} $\Lambda (w,\infty )\in {\cal
LI} (\Lambda (X), M^\ast (X); L^\infty )$
 if and only if
$$
\int_0^t{d\phi_X(x)\over w(x)}\leq
C{\phi_X(t)\over w(t)}, \qquad \forall t>0.\eqno { (2.3)}
$$

\P  First of all we remark that condition $(2.2)$ implies that $w$ satisfies
$(AM_q).$

{\bf (i)} We only have to prove
that the operator $Q_{\Lambda (X)}$ defined by $$Q_{\Lambda (X)}f=H(f \circ
\phi_X^{-1})\circ \phi_X $$
is bounded on $\Lambda (w,q)$, for nonnegative and nonincreasing functions.
$$\eqalign{
\V Q_{\Lambda (X)}f\V_{\Lambda (w,q)}^q
&=\int_0^\infty H(f \circ \phi_X^{-1})( \phi_X(t))^q w(t)dt \cr
&= \int_0^\infty H(f \circ \phi_X^{-1})( y)^q
w(\phi_X^{-1}(y))(\phi_X^{-1})^{\prime}(y)dy      \cr
& \leq C\int_0^\infty f ( \phi_X^{-1}(y))^q
w(\phi_X^{-1}(y))(\phi_X^{-1})^{\prime}(y)dy \cr
&= C
\int_0^\infty f(x)^qw(x)dx= C \V f\V_{\Lambda (w,q)}^q
}$$
where, by using  condition $(AM_q)$, the inequality  is satisfied if
and only if the weight
$$
v(y)=w(\phi_X^{-1}(y))(\phi_X^{-1})^{\prime}(y)
$$
satisfies
$$
\int_t^\infty{v(y)\over y^q}dx\leq{B\over t^q}\int_0^tv(y)dy \qquad \forall
t>0
$$
for some constant $B>0$ and this inequality is equivalent to $(2.2).$
\medskip
\noindent {\bf (ii)} This proof  is simpler.  Suppose that (2.3) holds. If
$f\in \Lambda (w,\infty )$ we have
$$\eqalign{
\sup_{t>0}{w(t)\over\phi_X(t)} \int_0^tf^\ast (x)d\phi_X(x)
&\leq
\V f\V_{\Lambda (w,\infty)}\
\sup_{t>0}{w(t)\over\phi_X(t)} \int_0^t{d\phi_X(x)\over w(x)} \cr
&\leq
C\V f\V_{\Lambda (w,\infty)}
}$$
For the converse implication,  consider for each $t>0$ the function

\noindent $f^\ast (x)={1\over w(x)}\chi_{[0,t]}(x)$ $ \in \Lambda (w,\infty)$.
Then  the inequality
$$
{w(t)\over\phi_X(t)}\int_0^tf^\ast (x)d\phi_X(x)\leq C
\sup_{t>0}f^\ast (t)w(t), \qquad \forall t>0
$$
implies (2.3).
\fin
\medskip
 These results should be compared with those appearing in
{\bf [18]}. In his paper Sharpley deals with the case of interpolation between
$(\Lambda (X_1),M(X_1) )$ and $(\Lambda (X_2),M(X_2) )$ and characterizes when
 $\Lambda^\alpha (Y)$ or $M(Y) $ are  interpolated spaces. Observe that
$\Lambda^\alpha (Y)$ and $M(Y) $ are particular cases of $\Lambda (w,q)$.

Note that the operator $Q_{\Lambda (X)}$ given by (2.1) is actually \lq\lq
linear". Hence, proposition 3 is more or less immediate in the sense  that
proposition 2 is not needed in order to prove proposition 3 since the
monotonicity condition supplies the linear operator $Q_{\Lambda (X)}$. That is,
the linearization made in proposition 2 is not needed.

However, our general result also allows us to treat the less evident case of
$(X,M(X))$. First of all, we will translate theorem 1 to this context:
\medskip
\proclaim Proposition 5. If $X$ is a r.i. Banach space then ${\cal
LI}(X, M^\ast (X); L^\infty )={\cal
QLI}(X,\break M^\ast (X); L^\infty )= {\cal
M}(X,M^\ast (X))={\cal
LPI}(X,M^\ast (X);L^\infty )$ (the last fact if $I=[0,1]$). Furthermore,
a couple
of spaces  $(Y,Z)$ belongs to any of the classes stated before  if and only if
the quasilinear operator $Q_X$ defined by
$$
Q_X\varphi (t)= { 1\over {\phi_X(t)}} \V \varphi\chi_{[0,t]}\V_X
$$
is bounded from $Y$ into $Z$ for nonnegative and nonincreasing functions.

\medskip

It is easy to prove that the quasilinear operator $Q_X$ satisfies
$$
Q_X\va (t)\leq {1\over \phi_X(t)}K(\phi_X(t),\va ;X,L^\infty )\leq 2Q_X\va (t)
$$
for all $t>0$ and for all nonnegative nonincreasing function $\va \in
X+L^\infty$, where $K$ is the $K$-functional introduced by Peetre.

 When $X=L^p$,
$1<p<\infty$, it is not difficult  to prove that for a $p$-convex space
 $Y$:   $(Y,Y)\in {\cal
LI}(L^p,L^{p,\infty}; L^\infty )$ $\Longleftrightarrow$
$(Y,Y)\in {\cal
LI}(L^{p,1},L^{p,\infty}; L^\infty )$,
 (see {\bf [9]} for a definiton of $p$-convex spaces).

In the next result we characterize the spaces $\Lambda (w,q)$ which are
interpolated in this context.
\medskip
\proclaim Proposition 6. Let $X$ be a r.i. Banach space and suppose that
$\Lambda
(w,q)$ is a Banach space, $1\leq q\leq \infty.$ Then the following is true:
\item{\bf (i)} If $1\leq q<\infty$, $\Lambda (w,q)\in {\cal
LI} (X, M^\ast (X); L^\infty )$ $\Longleftrightarrow$  condition (2.2) holds.
 \item{\bf (ii)} $\Lambda (w,\infty )\in {\cal
LI} (X, M^\ast (X); L^\infty )$
 if and only if the following condition  is satisfied
$$
{w(t)\over \phi_X(t)} \left\V {\chi_{[0,t]}\over w}\right\V_X\leq C  \qquad
\forall t>0.\eqno (2.4) $$

\P {\bf (i)}  We know that $\Lambda (w,q)\in {\cal
LI} (X, M^\ast (X); L^\infty )$ $\Longleftrightarrow$ the operator
$$
Q_X(f)(t)={1\over \phi_X(t)}\V f\chi_{[0,t]}\V_X
$$
is bounded in $\Lambda (w,q)$ for nonincreasing and nonnegative
functions. Since $Q_X(f)\leq Q_{\Lambda (X)}(f),$ we obtain that condition
(2.2) implies the interpolation property. On the other hand, since
$$
Q_X(\chi_{[0,s]})(t)=\chi_{[0,s]}(t)+{\phi_X(s)\over \phi_X(t)
}\chi_{[s,\infty)}
$$
we have that if $Q_X$ is bounded on $\Lambda (w,q)$, then
$$
\int_0^s w(x)dx + \phi_X(s)^q\int_s^\infty {w(x)\over\phi_X(x)^q}dx\leq
C \int_0^s w(x)dx
$$
and hence (2.2) holds.

\noindent {\bf (ii)} Suppose $\Lambda (w,\infty )\in {\cal
LI} (X, M^\ast (X); L^\infty )$. By observing that ${1\over w}\in \Lambda
(w,\infty
)$ and that
 $$
Q_X\left({1\over w}\right)(t)=
{1\over \phi_X(t)}\left\V{\chi_{[0,t]}\over w}\right\V_X
$$
we get condition (2.4).

For the reverse part we know that
 $w(x)f^\ast(x)\leq  \V f\V_{\Lambda (w,\infty )}$, $\forall x>0$ and $\forall
f\in \Lambda (w,\infty )$. Furthermore
 $$\eqalign{
Q_X(f^\ast)(t)&={1\over\phi_X(t)}\V f^\ast\chi_{[0,t]}\V_X \cr
\leq &{C\over\phi_X(t)}\V f\V_{\Lambda (w,\infty )}\left\V{\chi_{[0,t]}\over
w}\right\V_X \cr
& \leq {C\over w(t)}\V f\V_{\Lambda (w,\infty )}
 }$$
and since $w$ is nondecreasing we obtain that the operator $Q_X$ is bounded on
$\Lambda (w,\infty )$.
\fin
\medskip
{\noindent \bf Remarks.}
\noindent {\bf (i)} We observe from Propositions 4 and 6 that the spaces
$\Lambda (w,q)$, $ 1\leq q <\infty$, belong to ${\cal LI} (X, M^\ast (X);
L^\infty )$ and ${\cal LI} (\Lambda (X), M^\ast (X); L^\infty ) $
simultaneously. Furthermore, if $X=L^p$, $1<p<q$, the weight $w$
satisfies condition $(AM_{{q\over p})}$ and then the space $\Lambda (w,q)$ is
$p$-convex. The converse is also true, i.e., let $q>p\geq 1$, if the space
$\Lambda (w,q)$ is $p$-convex then $\Lambda (w,q)\in {\cal LI} (L^p,
L^{p,\infty}; L^\infty )$. This result should be compared with those appearing
in {\bf [16]}.

\noindent {\bf (ii)} If $\Lambda (X)\not= X=M^\ast (X)$ (for instance,
$X=L^{p,\infty},$ $p>1$) the space $X=\Lambda (\phi_X, \infty )\in {\cal
LI} (X, M^\ast (X); L^\infty )
$ but $X\notin {\cal
LI} (\Lambda(X), M^\ast (X); L^\infty )
$ as it is easy to prove.

\noindent {\bf iii)}   $\Lambda (w,\infty )\in {\cal LI}(L^{p,1}, L^{p,\infty};
L^\infty )$  $\Longleftrightarrow$ $\Lambda (w,\infty )$ is $p$-convex.
\medskip

In the case $X=L^p, 1<p<\infty$,   the situation is clearer  as the
following result shows
\medskip
\proclaim Proposition 7. Let $1\leq p<\infty$, a real
number and let $Y$ be a r.i. space. The following assertions are equivalent:
\item{\bf (i)} $Y\in {\cal LI}(L^{p,1}, L^{p,\infty};
L^\infty )$
\item{\bf (ii)} $Y\in {\cal LI}(L^{p,r}, L^{p,\infty}; L^\infty
)$, for some $1<r<\infty.$
\item{\bf (iii)} The upper Boyd index $\alpha (Y)< {1\over p}.$
(See {\bf [4]} for the definition of Boyd indices).

\P We only have to prove
$(ii)\Longrightarrow (i)$. We  restrict ourselves to the case $r=p$ because
the proof is similar. By proposition 6, we know that the operator
$$
Q_pf(t) = {1\over t^{1/p}}\left(\int_0^tf^p\right)^{1/p}
$$
is bounded in $Y$ for nonincreasing nonnegative functions, i.e. $\V
Q_pf\V_Y\leq C\V f\V_Y$, for some constant $C>0$. Let
 $f=f^{\ast}\in Y.$ It is quite easy
to compute that
$$
Q_p^{(n)}f(t) = \left(\int_0^1f^p(tx){[\log (1/x)]^n\over
n!}dx\right)^{1/p} $$
for any natural number $n\in {\bf N}$. If $\epsilon C<1$
we define
 $$
Sf(t)= \left( \sum_{n=0}^\infty [\epsilon^n
Q_p^{(n+1)}f(t)]^p\right)^{1/p}. $$
Since
$
S_Nf(t)= \left( \sum_{n=0}^N [\epsilon^n
Q_p^{(n+1)}f(t)]^p\right)^{1/p} \in Y
$, $S_Nf(t)\uparrow Sf(t)$ when $N\rightarrow\infty$, and
$$\eqalign{
\V S_Nf\V_Y&\leq \left\V \sum_0^N\epsilon^nQ_p^{(n+1)}f\right\V_Y    \cr
&\leq \left( \sum_0^\infty \epsilon^nC^{n+1}\right)  \V f\V_Y
}$$
we obtain that $Sf\in Y$ and $\V Sf\V_Y\leq C'\V f\V_Y$. But
$$\eqalign{
Sf(t)= & \left( \sum_0^\infty \int
_0^1f(tx)^p{[\epsilon\log(1/x)]^n\over n!}dx\right)^{1/p}  \cr
&=   \left(\int
_0^1f(tx)^p\sum_0^\infty{[\epsilon\log(1/x)]^n\over
 n!}dx\right)^{1/p} \cr
&= \left( \int_0^1f(tx)^p{dx\over x^\epsilon} \right)^{1\over p} =
{1\over t^{1-\epsilon\over p}}\left(
\int_0^tf(x)^pd(x^{1-\epsilon\over p})
\right)^{1/p}
 \cr
&= {1\over t^{1-\epsilon\over p}} \left\V f\chi_{[0,t]}\right\V_{L^
{{p\over 1-\epsilon},p}}
}$$
Hence $(Y,Y)\in {\cal LI}(L^{{p\over 1-\epsilon},p}, L^{p,\infty};
L^\infty )$ and so (i) is true.
\fin
\medskip

In the last part of this paper we consider a similar situation to the previous
one in the framework of Lorentz-Orlicz spaces. Different versions of this
class of spaces appear in {\bf [12]}, {\bf [19]} and they have been also
studied in {\bf [7]}, {\bf [13]}, {\bf [14]} and {\bf [16]}.   Here we
 consider
the Lorentz-Orlicz spaces as they appear in {\bf [12]}.

In the sequel $\va$ will denote an Orlicz  function, i.e. a convex,
non-decreasing function on $[0,\infty )$ such that $\va (0)=0$ and
$\lim_{t\to \infty}\va (t)=\infty.$ We also suppose that $\va$ satisfies the
$\Delta_2$ condition: there exists a constant $C>0$  such that $\va (2t)\leq
C\va (t),$ for all $t>0$, or equivalently,   there
exists $1\leq q<\infty$,  so that,
$$
\va (at)\leq a^q \va (t),\qquad \forall a\geq 1, t>0.
$$

The weight $w$ is supposed to satisfy the same conditions appearing in the
definition of Lorentz spaces, namely, $w$ is an a.e.
positive weight defined on $[0,\infty )$ such that $\int_0^tw<\infty,
\forall t<\infty $ and $\int_0^\infty w=\infty$.

The space $\Lambda (w, \va )$ is  the class of real valued
measurable functions on $I$ so that
$
\int_I\va (f^\ast (t))w(t)dt<\infty .
$

Next results are an extension of Ari\~no-Muckenhoupt's inequality.

\medskip

\proclaim Lemma 3. Let $\va$ be an Orlicz function. Suppose that there
exists a
constant $B>0$ such that
 $$
\int_t^\infty \va\left( {at\over x}\right)w(x)dx
\leq B \va (a)\int_0^tw(x)dx\qquad
\forall t>0, \forall a>0 \eqno (A_\varphi ) $$
then, for some $\alpha <1$ and  $D>0,$ we have
$$
\int_t^\infty \psi\left( {at\over x}\right)w(x)dx\leq D \psi (a)
\int_0^tw(x)dx
\qquad \forall t>0, \forall a>0
$$
where $\psi (t)=\va (t^\alpha ).$

\P We can
adapt
the arguments in {\bf [1]} to our more general situation. Only a few changes
are
necessary. Using the notation in  {\bf [1]} the number $\alpha <1 $ is choosen
in such a way that
 $$
2S^{1/q}=2\left({2B+1\over 2B+2}\right)^{1/q}< 2^\alpha <2.
$$
\hfill \fin
\medskip

As a consequence of this lemma and  with the same reasons appearing  in {\bf
[1]},
 we
obtain the following

\medskip

\proclaim Proposition 8.   Let $\va$ be an Orlicz
function. A weight $w$ satisfies condition $(A_\varphi )$ if and
only if there exists a constant $B'>0$ such that for every nonnegative
nonincreasing function $f$ on $(0,\infty )$ we have
$$
\int_0^\infty \va (H(f)) w\leq B' \int_0^\infty \va (f)w
$$
where {H(f)} is the Hardy transform of $f$.

\medskip

{\bf \noindent Remark.} If the weight satisfies $(A_\varphi )$,
 the expression
 $$
\V f\V_{\Lambda (w,\va )} = \inf \left\{ \rho ; \int_0^\infty \va\left({f^\ast
\over \rho}\right) w\leq 1 \right\}
$$
defines a quasi-norm on $\Lambda (w,\va )$ which is equivalent to the norm
$$|||f|||=\V Hf^\ast\V_{\Lambda (w,\va )}$$
 and
therefore  $\Lambda (w,\va )$ is a Banach space.

By introducing the Simonenko indices (see {\bf [12]}) we can
give neccesary or sufficent conditions in order to verify
condition $(A_\varphi )$ for a weight. Given an Orlicz convex function $\va$
 and
a number $T> 0$ we define
 $$p_T=\inf_{t\geq T} {t\va^{\prime}(t)\over \va (t)}$$
$$
q_T = \sup_{t\geq T} {t\va^{\prime}(t)\over \va (t)}
$$
where $\va^{\prime}(t)$ is supposed to be the right derivative of the Orlicz
function $\va$. We also introduce $p_0= \inf p_T$ and $q_0=\sup q_T.$ It is
clear that $1\leq p_0\leq p_T\leq q_T \leq q_0<\infty$ and
$$
\alpha^{q_T}\va
(t)\leq \va (\alpha t)\leq \alpha^{p_T}\va (t)$$
whenever $\alpha \leq 1$ and $T\leq \alpha t$.
\medskip
\noindent {\bf Proposition 9.} {\sl Let $\va $ be an Orlicz function and let
 $w$
be a weight. The following assertions are true:

{\bf i)}  If  a weight $w$ satisfies the condition $(AM_{p_0})$
then it also verifies $(A_\va)$. Furthermore
 the Hardy transform is bounded on
$\Lambda (w, \va )$ and $\Lambda (w, \va )$ is a Banach space.

{\bf ii)} If the Hardy transform is bounded on
$\Lambda (w, \va )$ then $w$ has to satisfy condition $(AM_{q_0})$}.

\medskip

\P {\bf i)}  First of all suppose that $w$ satisfies the
condition $(AM_{p_0})$, then
$$
\int_t^\infty \va \left( {at\over x}\right) w(x)dx\leq \va (a)
\int_t^\infty  \left({t\over x}\right)^{p_0}w(x)dx $$
and hence $w$ verifies condition $(A_\va )$ and so
the Hardy transform is bounded on
$\Lambda (w, \va )$. The other statements are clear.

{\bf ii)} Now we assume that
 there exists a constant
$C\geq 1$ such that
$$
\V Hf\V_{\Lambda (w,\va )}\leq C \V f\V_{\Lambda (w,\va )}
$$
for all non increasing function $f\in \Lambda (w, \va )$. In particular,
 given
$t>0$ if $s=\V \chi_{[0,t]}\V_{\Lambda (w,\va )}^{-1}$, we obtain that
$$
\V H(\chi_{[0,t]})\V_{\Lambda (w,\va )}\leq {C\over s}.
$$
Therefore
$$
\int_0^t\va\left({s\over C}\right)w(x)dx+
\int_t^\infty\va\left({st\over Cx}\right)w(x)dx\leq 1.
$$
Since $\displaystyle {\va (s)\over C^{q_0}}\leq \va\left({s\over C}\right)$
 and
$$
 \va\left({st\over Cx}\right)\geq \va\left({s\over C}\right)
\left({t\over x}\right)^{q_0}\geq \va (s)
\left({t\over C x}\right)^{q_0}
$$
we have that
$$
\int_0^tw(x)dx+\int_t^\infty\left({t\over x}\right)^{q_0}w(x)dx
\leq C^{q_0}  \int_0^tw(x)dx.
$$
and this completes the proof.
\fin
\medskip

In the case of considering the space $\Lambda (w,\va )$ on the unit interval
$I=[0,1]$ (or more generally $I=[0,l]$ for $l<\infty$) we can give
a more precise
result.
\medskip
\noindent{\bf Proposition 10}. {\sl Suppose $I=[0,l]$, ($l<\infty$) and let
 $\va$
an Orlicz function. Let
$ p=\liminf_{t\to\infty} {t\va^{\prime}(t)\over
\va (t)}$  and
 $q=\limsup_{t\to\infty}
{t\va^{\prime}(t)\over \va (t)}.$
The following assertions are true:

{\bf i)} If $w$ satisfies condition $(AM_p)$ then the Hardy transform is
bounded on $\Lambda (w,\va )$ and $\Lambda (w,\va )$ is a Banach space.

{\bf ii)} If the Hardy transform is
bounded on $\Lambda (w,\va )$ then $w$ satisfies condition $(AM_{q+\epsilon})$
for all $\epsilon >0.$ }

\medskip
\P {\bf i)} If $w$
 satisfies condition $(AM_p)$ by using Lemma 2.1 in {\bf [1]} $w$ also verifies
condition  $(AM_{p-\epsilon })$ for some $\epsilon >0.$ Then there exists $T>0$
such that $p-\epsilon <p_T.$ Define now the function $\var_T$ by
$$
\var_T (t)=\cases
{
\va (t), & if $t\geq T$ ; \cr
\va (T)\left({t\over T}\right)^{p_T}, & otherwise. \cr
}
$$
Note that the function $\var_T $ is an Orlicz function
for which
$$
p_T=\inf_{t>0}{t\var^{\prime}(t)\over \var (t)}.
$$
Hence by using preceding proposition  we obtain that Hardy transform is
bounded
on
$\Lambda (w,\var_T ).$ Thus  the Hardy transform is
also bounded on $\Lambda (w,\va )$,
 as  $m(I)<\infty$ and the Orlicz functions are equivalent at infinity.

{\bf ii)} Let $\epsilon $ be a positive number. There exists $T>0$ such that
$$
{t\va^{\prime}(t)\over \va (t)} \leq q_T< q+\epsilon \qquad \forall t\geq T
$$
If we define the Orlicz function
$\var_T$ as before it is clear that $\sup_{t>0}{t\var^{\prime}(t)\over \var
(t)}\leq q_T$. Since $\var_T$ is equivalent to $\va$ at infinity  and the
Hardy transform is bounded on $\Lambda (w, \var_T )$, by using again
Proposition 9, we deduce that $w$ satisfies the  condition $(AM_{q_T} )$ and
so
$(AM_{q+\epsilon} )$.

This completes the proof. \fin
 \medskip

  Now we can state the
corresponding interpolation results  whose proves follow the same lines as in
Propositions 4 and 6.
 \medskip
\noindent {\bf Proposition 11}. {\sl Suppose that $\va$ is an Orlicz function
and $w$ satisfies the condition  $(A_\va
)$. Let $X$ be a r.i. function space. Then the following assertions are
equivalent:}

{\sl \noindent{\bf i)} Any linear or cuasilinear operator $T$ which is bounded
from $\Lambda (X)$ into $M^\ast (X)$ and from $L^\infty$ into $L^\infty$
verifies
$$
\int _0^\infty \va ((Tf)^\ast(x))w(x)dx \leq C'
\int _0^\infty \va (f^\ast (x))w(x)dx
$$
for some constant $C'>0$ and for any function $f\in \Lambda (w, \va).$

\noindent {\bf ii)} The same as in i) but with operators
mapping $X$ into $M^\ast
(X)$, instead of
$\Lambda (X)$ into $M^\ast (X)$.

\noindent {\bf iii)} There exist a constant $D>0$ such that
$$
\int_t^\infty \va\left({\phi_X(t)\over \phi_X(x)}\right)w(x)dx
\leq D \int_0^t w(x)dx
$$

\noindent Moreover,  if one of these conditions is satisfied for $X=L^p$,
$1\leq p< \infty$, the space $\Lambda (w,\va )$ is $p$-convex.}

 \bigskip
{\noindent \bf Acknowledgment}

We want to thank  Stephen Montgomery-Smith for polishing an earlier version of
the condition $(A_\va )$.
\bigskip
{\noindent \bf References}
\medskip
\item{\bf [1]}  M. A. Ari\~no and B. Muckenhoupt, {\it Maximal Functions on
Classical Lorentz Spaces and Hardy's Inequality with Weights for Nonincreasing
Functions}. Trans. A.M.S. {\bf 320}, (2), (1990), 727-735.

\item{\bf [2]} C. Bennet and R. Sharpley, Interpolation of operators. Academic
Press, 1988.

\item{\bf [3]} J. Bergh and J. L\"ofstr\"om, Interpolation Spaces. An
Introduction. Springer Verlag,  1976.

\item{\bf [4]} D. W. Boyd, {\it Indices of function spaces and their
 relationship
to interpolation}. Canad. J. Math. {\bf 21}  (1969), 1245-1254.

\item{\bf [5]} A. P. Calder\'on,  {\it Spaces between $L^1$ and $L^\infty$ and
the theorem of Marcinkiewicz}. Studia Math. {\bf  26}  (1966), 273-299.

\item{\bf [6]} J. Garcia Cuerva and J. L. Rubio de Francia, Weighted norm
inequalities and related topics. North-Holland Mathematical Studies, Amsterdam,
1985

\item{\bf [7]} N. J.  Kalton,  {\it Endomorphisms of Symmetric Function
 Spaces}.
Indiana Univ. Math. J. {\bf 34}. No. 2 (1985), 225-247.

\item{\bf [8]} A. Kamminska. {\it Some remarks on Orlicz-Lorentz spaces}.
Math. Nach. {\bf 147}, (1990), 29-38.

\item {\bf [9]} J. Lindenstrauss and L. Tzafriri, Classical Banach Spaces II.
Springer Verlag, 1979.

\item{\bf [10]} G. Lorentz and T. Shimogaki,  {\it Interpolation Theorems for
the
Pairs of Spaces $(L^p, L^\infty )$ and $(L^1, L^q)$ }. Trans. A.M.S. {\bf 139}
  (1971), 207-221.

\item{\bf [11]} L. Maligranda, {\it A generalization of the Shimogaki theorem}.
Studia Math. {\bf 71}  (1981), 69-83.

\item{\bf [12]} L. Maligranda, {\it Indices and interpolation}. Dissertationes
Math. {\bf 234}  (1985), 1-54.

\item{\bf [13]} M. Mastylo, {\it Interpolation of linear operators in
Calder\'on-Lozanovskii spaces}.
Comment. Math. {\bf 26}, 2, (1986), 247-256.

\item{\bf [14]} S. J. Montgomery-Smith, {\it Comparison of Orlicz-Lorentz
spaces}. Preprint.

\item{\bf [15]} B. Muckenhoupt, {\it Weighted norm inequalities for the Hardy
maximal function}. Trans. Amer. Math. Soc. {\bf 165}, (1972), 207-226.

\item{\bf [16]} Y. Raynaud, {\it On Lorentz-Sharpley
spaces}. Proceedings of the Workshop \lq\lq Interpolation Spaces and Related
Topics", Haifa, June 1990.

\item {\bf [17]} E. Sawyer {\it  Boundedness of Classical Operators on
Classical
Lorentz Spaces}. Studia Math. {\bf 96}, (1990), 145-158.

\item{\bf [18]} R. Sharpley, {\it Spaces $\Lambda_\alpha (X)$ and
Interpolation}.
Journal of Functional Analysis {\bf 11},   (1972), 479-513.

\item{\bf [19]} A. Torchinsky, {\it Interpolation of operators and Orlicz
classes}. Studia Math. {\bf 59}, (1976), 177-207.

\item{\bf [20]} T. Shimogaki {\it An interpolation theorem on Banach function
spaces}. Studia Math. {\bf 31},  (1968),  233-240.

\item{\bf [21]} M. Zippin, {\it Interpolation of Operators of  Weak Type
Between
Rearrangement Invariant Function Spaces}. Journal of Functional Analysis
{\bf 7
},   (1971), 267-284.

\bigskip

Jes\'us Bastero and Francisco J. Ruiz

{\fp Departamento de Matem\'aticas

Facultad de Ciencias

Universidad de Zaragoza

50009 Zaragoza (Spain)}

\bye